\documentclass{amsart}
\usepackage{amsfonts}
\usepackage{amsmath,amssymb}
\usepackage{amsthm}
\usepackage{amscd}
\usepackage{graphics}
\usepackage{graphicx}
\theoremstyle{definition}{
	\newtheorem{Def}{{\rm Definition}}
	\newtheorem{Ex}{{\rm Example}}
	\newtheorem{Rem}{{\rm Remark}}
	\newtheorem{Prob}{{\rm Problem}}
}
\theoremstyle{plain}
{
	\newtheorem{Cor}{Corollary}
	\newtheorem{Prop}{Proposition}
	\newtheorem{Thm}{Theorem}
	\newtheorem{MainThm}{Main Theorem}

}

\begin{document}
	\title[Special generic maps can present nice generalized multisections]{Manifolds admitting special generic maps and their nice generalized multisections}
	\author{Naoki Kitazawa}
	\address{Institute of Mathematics for Industry, Kyushu University, 744 Motooka, Nishi-ku Fukuoka 819-0395, Japan\\
		TEL (Office): +81-92-802-4402 \\
		FAX (Office): +81-92-802-4405 \\
	}
	\email{n-kitazawa@imi.kyushu-u.ac.jp}
	\urladdr{https://naokikitazawa.github.io/NaokiKitazawa.html}
	
	\begin{abstract}
		We show that manifolds admitting {\it special generic} maps also admit nice {\it generalized multisections}.
		
		{\it Special generic} maps are natural generalized versions of Morse functions with exactly two singular points on closed manifolds, characterizing spheres whose dimensions are not $4$ topologically and the $4$-dimensional unit sphere, and canonical projections of unit spheres. They are shown to restrict the differentiable structures of spheres etc. and topologies of more general manifolds strongly by Saeki, Sakuma etc., followed by Nishioka, Wrazidlo etc. and followed by the author. Some elementary or important manifolds also admit such maps. 
		
		({\it Generalized}){\it multisections} of manifolds are nice decompositions of (compact) manifolds, generalizing so-called Heegaard splittings of $3$-dimensional manifolds. 
		PL manifolds have been shown to have (generalized) multisections enjoying certain properties by Rubinstein and Tillmann.


		
	\end{abstract}
	
	
	\maketitle
	\section{Introduction.}
	\label{sec:1}
	A Morse function with exactly two {\it singular} points on a closed manifold is one of simplest maps in Morse theory and differential topology of manifolds. This characterize spheres whose dimensions are not $4$ and the $4$-dimensional unit sphere. This is so-called Reeb's theorem. A {\it special} generic map is a higher dimensional generalization and a central object in our paper. Canonical projections of unit spheres are special generic.
	
	Decomposing manifolds is a fundamental, natural and strong tool in understanding their (global) geometric properties.
	A ({\it generalized}) {\it multisection} of a compact manifold is a kind of nice decompositions of a manifold. A {\it Heegaard} splitting of a $3$-dimensional compact and connected manifold is a specific case. This presents a big and interesting topic on $3$-dimensional manifolds. See \cite{hempel} for fundamental theory on these manfiolds. $4$-dimensional versions are well-known. \cite{gaykirby} is a pioneering study. \cite{rubinsteintillmann1, rubinsteintillmann2} presents one of generalized cases for manifolds of general dimensions.
	
	We introduce fundamental notions, terminologies and notation which are important in our paper and present our Main Theorem. 
	
	For a positive integer $k$, ${\mathbb{R}}^k$ denotes the $k$-dimensional Euclidean space, which is a $k$-dimensional smooth manifold of course and a Riemannian manifold with the standard Euclidean metric. $||x|| \geq 0$ denotes the distance betwen the origin $0 \in {\mathbb{R}}^k$ and $x \in {\mathbb{R}}^k$. $S^k:=\{x \in {\mathbb{R}}^{k+1} \mid ||x||=1 \}$ denotes the $k$-dimensional unit sphere for $k \geq 0$. It is a $k$-dimensional smooth closed submanifold of ${\mathbb{R}}^{k+1}$ with no boundary. $D^k:=\{x \in {\mathbb{R}}^{k} \mid ||x|| \leq 1 \}$ denotes the $k$-dimensional unit disk for $k \geq 1$. It is a $k$-dimensional smooth closed submanifold of ${\mathbb{R}}^{k+1}$ whose boundary is $S^{k-1}$.

	We omit rigorous expositions on fundamental terminologies, notions and notation on complexes, polyhedra, and other related important ones related to our paper.

	For systematic geometric or combinatorial theory on the PL category and the piesewise category, which is equivalent to the PL category, see \cite{hudson} for example.
	
	For a topological space $X$ homeomorphic to a cell complex whose maximal dimension is finite, $\dim X$ denotes its dimension, which is uniquely defined.
	
	Manifolds are homeomorphic to CW complexes. Every smooth manifold is canonically regarded as a so-called {\it PL manifold}.
	
	For a differentiable map $c:X \rightarrow Y$ between differentiable manifolds, $p \in X$ is a {\it singular} point of $c$ if the rank of the differential $dc_p$ at $p$ is smaller than both $\dim X$ and $\dim Y$. $c(p)$ is also a {\it singular value} of $c$. $S(c)$ denotes the set of all singular points of $c$ and we call this the {\it singular set} of $c$.
	
	A {\it diffeomorphism} between smooth manifolds means a smooth map with no singular points which is also a homeomorphism. A {\it diffeomorphism on a smooth manifold} means a diffeomorphism from the manifold onto itself. Two manifolds are defined to be {\it diffeomorphic} if and only if there exists a diffeomorphism between these two. This of course defines an equivalence relation on the class of all smooth manifolds (where corners are eliminated). We can also define {\it PL homeomorphisms} and {\it PL homeomorphic manifolds} via piesewise smooth homeomorphisms similarly. 
	
	\begin{Def}
		\label{def:1}
		A smooth map $c:X \rightarrow Y$ between two smooth mnaifolds with no boudnaries are said to be {\it special generic} if at each singular point $p \in X$ there exists suitable local coordinates around $p$ and $f$ is represented as $(x_1,\cdots,x_{\dim X}) \rightarrow (x_1,\cdots,x_{\dim Y-1},{\Sigma}_{j=1}^{\dim X-\dim Y+1} {x_{\dim Y+j-1}}^2)$.
	\end{Def}
	Later, Proposition \ref{prop:1} says that the image of a special geneirc map on a closed manifold (into a Euclidean space) is a smoothly immersed manifold of codimension $0$. 
	It is a kind of fundamental exercises on smooth manifolds, Morse functions, and singularity theory of smooth maps, to see that the canonical projection of a unit sphere, defined as a map mapping $(x_1,x_2) \in S^k \subset {\mathbb{R}}^{k+1}={\mathbb{R}}^{k_1} \times {\mathbb{R}}^{k_2}$ into $x_1 \in {\mathbb{R}}^{k_1}$ where $k \geq 2$, $k_1,k_2 \geq 1$ and $k=k_1+k_2$. As pioneering stuides, \cite{burletderham, furuyaporto} are known. 
	Since 1990s, manifolds admitting special generic maps have been studied by Saeki and Sakuma as \cite{saeki1,saeki2,saekisakuma1,saekisakuma2}, followed by \cite{nishioka,wrazidlo1,wrazidlo2,wrazidlo3}.
	They have revealed restrictions on the differentiable structures of the homotopy spheres etc. and the homology groups. As a pioneer, the author has studied the cohomology rings of the manifolds mainly in \cite{kitazawa1,kitazawa2,kitazawa3,kitazawa4,kitazawa5,kitazawa6,kitazawa7,kitazawa8} for example.
	
	Let $k$ be a positive integer and $m$ another integer greater than $1$.
	A ({\it generalized}) {\it multisections} of a smooth or PL compact manifold of {\it degree $k$} is a decomposition of an $m$-dimensional compact manifold into $k$ {\it $1$-handlebodies} (resp. compact and connected manifolds collapsing to lower dimensional polyhedra). A {\it $1$-handlebody} means a manifold diffeomorphic to $D^m$ or represented as a boundary connected sum of finitely many copies of $S^1 \times D^{m-1}$ (considered in the smooth category). This decomposition can be discussed in the PL category and the smooth category. We also consider the classes of {\it normal} decompositions for these decompositions for example.
	Some of these expositions are due to  \cite{rubinsteintillmann1, rubinsteintillmann2} whereas some are arranged according to our own thought on these notions.
	We introduce a specific class of generalized multisections generalizing the class of multisections in a natural way later as the class of {\it near multisections}. This is also a new ingredient of our paper.  
	
For example,  \cite{rubinsteintillmann1, rubinsteintillmann2} show that every $2k$-dimensional or ($2k+1$)-dimensional PL closed and connected manifold admits a PL {\it neat} multisection of degree $k+1$ where $k$ is an integer satisfying $2k>0$ or $2k+1>0$ respectively. Our Main Theorems give new answers on existence and explicit examples of such decompositions as follows, giving answers to our Problems \ref{prob:1} and \ref{prob:2}, proposed later. A smooth ot PL ({\it near}) {\it neat} multisection is defined in Definition \ref{def:4}. 
	\begin{MainThm}
		\label{mthm:1}
		Let an $m$-dimensional closed, connected and orientable manifold $M$ admit a special generic map $f:M \rightarrow {\mathbb{R}}^n$ with $m>n$ such that the image is a smoothly immersed $n$-dimensional compact and connected manifold{\rm :} let this $n$-dimensional smooth manifold be obtained by removing the interiors of smoothly and disjointly embedded copies of the unit disk $D^{n}$ from a smooth closed one admitting a smooth {\rm (}PL{\rm )} near neat multisection of degree $k$. Then $M$ admits a smooth {\rm (}resp. PL{\rm )} near multisection of degree k.
		In the case the given multisection is neat, the resulting near multisection is normal.
	\end{MainThm}
Main Theorem \ref{mthm:2}, presented in the fourth section, is an application of Main Theorem \ref{mthm:1}. Our work is also strongly motivated by the following present situations.
	\begin{itemize}
		\item Before the study and other than \cite{rubinsteintillmann1, rubinsteintillmann2}, existence of (generalized) multisections under suitable situations and general construction have been known. However, it is a different and difficult study to find explicit cases of such structures. Examples and explicit or general theory of Heegaard splittings such as \cite{johnson1,johnson2} and that of multisections of $3$-dimensional manifolds such as \cite{gomezlarranaga,itoogawa,koenig,ogawa1,ogawa2} and $4$-dimensional variants such as \cite{koenig,meier} have been presented and related problems are still actively studied. Higher dimensional cases seem to have various new, interesting and related problems.
		\item \cite{baykursaeki1,baykursaeki2} present algorithms for obtaining so-called {\it generic} smooth maps into the plane from 4-dimensional compact manifolds respecting the structures of some multisections of degree $3$ or {\it trisections}. They are regarded as variants of studies related to understanding of Heegaard splittings via nice Morse functions. 
		In the pioneering work \cite{gaykirby}, deformations of Morse functions are important and related tools are important in studying such structures in the smooth category. \cite{asano,hayano} are also related recent studies.
		
		We try to understand such structures of nice classes via generic smooth maps of explicit and nice classes. Morse functions are regarded as simplest generic maps. Special generic maps are also regarded as very explicit maps.
		
		For generic smooth maps, see \cite{golubitskyguillemin} for example. 
		This explains about fundamental singularity theory of differentiable maps systematically and related some advanced topics.
		Note that {\it generic} smooth maps are not defined rigorously in general scenes. We may regard "generic" means "general in situations we discuss". So-called {\it stable} maps are regarded as generic for example, in the context of singularity theory of smooth maps. Consult also this book for stable maps. 
	\end{itemize}
	
	In the next section, we present additional fundamental properties of special generic maps. 
	In the third section, we introduce some classes of (generalized) multisections referring to \cite{rubinsteintillmann1, rubinsteintillmann2} or arranged ones. After that we introduce {\it near multisections} of ({\it degree k}) of manifolds. Defining this class is also an important ingredient of our new work.
	The fourth section is devoted to Main Theorems. \\
	\ \\
	{\bf Conflict of Interest.} \\
	The author is a member of the project JSPS KAKENHI Grant Number JP22K18267 "Visualizing twists in data through monodromy" (Principal Investigator: Osamu Saeki). The present study is due to this project. \\
	\ \\
	{\bf Data availability.} \\
	Data supporting our present study essentially are all in our paper.
	
	\section{Fundamental properties and existing studies on special generic maps and the manifolds.}

	The {\it diffeomorphism group} of a smooth manifold is the group of all diffeomorphisms on it whose topology is the so-called {\it Whitney $C^{\infty}$ topology}. \cite{golubitskyguillemin} explains about such topologies on spaces of smooth maps between given two smooth manifolds as a fundamental notion.
	
	A {\it smooth} bundle means a bundle and whose fiber is a smooth manifold whose structure group is the diffeomorphism group. A {\it linear} bundle means a bundle whose fiber is a Euclidean space, unit sphere, or a unit disk and whose structure group consists of linear transformations. Linear transformations are defined in a natural and canonical way.
	
	For general theory of bundles, see \cite{steenrod}. \cite{milnorstasheff} concentrates mainly on linear bundles or so-called vector bundles, which are specific linear bundles.
	\begin{Prop}
		\label{prop:1}
		A special generic map $f:M \rightarrow {\mathbb{R}}^n$ on an $m$-dimensional closed and connected manifold $M$ enjoys the following properties.
		\begin{enumerate}
			\item \label{prop:1.1}
			We have an $n$-dimensional compact and connected smooth manifold $W_f$, a smooth surjection $q_f:M \rightarrow W_f$, a smooth immersion $\bar{f}:W_f \rightarrow {\mathbb{R}}^n$ and the relation $f=\bar{f} \circ q_f$. $q_f$ maps the singular set $S(f)$ of $f$ onto the boundary $\partial W_f \subset W_f$ as a diffeomorphism.
			\item \label{prop:1.2}
			We have a small collar neighborhood $N(\partial W_f)$ of the boundary $\partial W_f \subset W_f$. Furthermore, we can have one such that the composition of the restriction of $q_f$ to the preimage ${q_f}^{-1}(N(\partial W_f))$ with the canonical projection to $\partial W_f$ gives a linear bundle whose fiber is the {\rm (}$m-n+1${\rm )}-dimensional unit disk $D^{m-n+1}$.
			\item \label{prop:1.3}
			The restriction of $q_f$ to the preimage of $W_f-{\rm Int}\ N(\partial W_f)$ gives a smooth bundle whose fiber is an {\rm (}$m-n${\rm )}-dimensional standard sphere. In some specific case, such as the case $m-n=0,1,2,3$, the bundle is linear.
		\end{enumerate}
	\end{Prop}
	\begin{Prop}
		\label{prop:2}
		Given a smooth immersion ${\bar{f}}_N:\bar{N} \rightarrow {\mathbb{R}}^n$ of an $n$-dimensional smooth, compact and connected manifold $\bar{N}$. We have a special generic map $f:M \rightarrow {\mathbb{R}}^n$ on a suitable $m$-dimensional closed and connected manifold $M$ enjoying the properties {\rm (}\ref{prop:1.1}{\rm )}, {\rm (}\ref{prop:1.2}{\rm )} and {\rm (}\ref{prop:1.3}{\rm )} of Proposition \ref{prop:1} and the following two where we abuse the notation.
		\begin{enumerate}
			\item \label{prop:2.1}
			The linear bundle of Proposition \ref{prop:1} {\rm (}\ref{prop:1.2}{\rm )} is a trivial linear bundle.
			\item \label{prop:2.2}
			The smooth bundle of Proposition \ref{prop:1} {\rm (}\ref{prop:1.3}{\rm )} is a trivial smooth bundle.
		\end{enumerate}
	\end{Prop}
	Other than canonical projections of unit spheres, we explain about simplest special generic maps.
	\begin{Ex}
		\label{ex:1}
		Let $l>0$ be an integer and $m \geq n \geq 2$ integers.
		We choose an integer $1 \leq n_j \leq n-1$ for each integer $1 \leq j \leq l$. We consider a connected sum of $l>0$ manifolds in the family $\{S^{n_j} \times S^{m-n_j}\}_{j=1}^l$ in the smooth category. We have a special generic map $f:M \rightarrow {\mathbb{R}}^n$ on the resulting manifold $M$ such that in Proposition \ref{prop:2}, the the map $\bar{f}$ is an embedding.
		
	\end{Ex}
	
	A {\it homotopy sphere} means a smooth manifold of dimensions at least $1$ which is homeomorphic to a unit sphere. 
	A {\it standard} (an {\it exotic}) sphere is a homotopy sphere which is diffeomorphic to some unit sphere (resp. not diffeomorphic to any unit sphere).
	
	\begin{Thm}[\cite{saeki1,saeki2}]
		\label{thm:1}
		\begin{enumerate}
			\item
			\label{thm:1.1}
			An $m$-dimensional closed and connected manifold $M$ of dimension $m \geq 2$ admits a special generic map into ${\mathbb{R}}^2$ if and only if either of the following holds.
			\begin{enumerate}
				\item $M$ is a homotopy sphere which is a homotopy sphere of dimension $m \neq 4$ or a standard sphere of dimension $m=4$.
				\item A manifold
				represented as a connected sum of smooth manifolds in the smooth category where each of the manifolds here is either of the following manifolds. 
				\begin{enumerate}
					\item The total space of a smooth bundle over $S^1$ whose fiber is a homotopy sphere and $m \neq 5$.
					\item The total space of a smooth bundle over $S^1$ whose fiber is a standard sphere and $m=5$.
				\end{enumerate}
			\end{enumerate}
			\item
			\label{thm:1.2}
			If an $m$-dimensional closed and simply-connected manifold $M$ of dimension $m \geq 4$ admits a special generic map into ${\mathbb{R}}^3$, then we have either of the following two.
			\begin{enumerate}
				\item $M$ is a homotopy sphere which is a homotopy sphere of dimension $m>4$ or a standard sphere of dimension $m=4$.
				\item $M$ is a manifold
				represented as a connected sum of smooth manifolds in the smooth category where each of the manifolds here is either of the following manifolds. 
				\begin{enumerate}
					\item The total space of a smooth bundle over $S^2$ whose fiber is a homotopy sphere and $m \neq 6$.
					\item The total space of a smooth bundle over $S^2$ whose fiber is a standard sphere and $m=6$.
				\end{enumerate}
			\end{enumerate}
			In the case $m=4,5$, the converse also holds where a fiber of each bundle is an {\rm (}$m-2${\rm)}-dimensional standard sphere.
			\item
			\label{thm:1.3}
			Furthermore, for any manifold in Theorem \ref{thm:1} {\rm (}\ref{thm:1.1}{\rm )} and {\rm (}\ref{thm:1.2}{\rm )}, we can construct a special generic map enjoying the properties of Example \ref{ex:1} except the triviality of the smooth bundle and the linear bundle in Proposition \ref{prop:1} {\rm (}\ref{prop:1.2}{\rm )} and {\rm (}\ref{prop:1.3}{\rm )} or Proposition \ref{prop:2} under the condition $n_j=1$ for any $j$ in {\rm (}\ref{thm:1.1}{\rm )} and under the condition $n_j=2$ for any $j$ in {\rm (}\ref{thm:1.2}{\rm )}.
		\end{enumerate}
	\end{Thm}

	\section{Multisections, generalized multisections and our new class of generalized multisections.}
	\begin{Def}
		\label{def:2}
		Let $k>1$ be an integer. A $k$-dimensional {\it $1$-handlebody} means a smooth manifold diffeomorphic to the $k$-dimensional unit disk $D^k$ or one represented as a boundary connected sum of finitely many copies of $S^1 \times D^{k-1}$ where the boundary connected sum is considered in the smooth category.
		
	\end{Def}
	Note again that a smooth manifold is regarded as a PL manifold canoncially. $1$-handlebodies are also PL manifolds of course.
	Here, 1-handlebodies and more general compact manifolds with non-empty boundaries may have suitable non-empty corners. Note that smooth manifolds with corners can be smoothed to unique smooth manifolds with no corners in canonical ways. 
	
	For 1-handlebodies, we can consider non-orientable ones. However we do not consider such cases here.

	The following definition is partially due to \cite{rubinsteintillmann1, rubinsteintillmann2} and some conditions and rules may be different from ones there.
	
	\begin{Def}
		\label{def:3}
		\begin{enumerate}
			\item 
			\label{def:3.1}
			For a smooth (PL), closed and connected manifold $X$, a family of $k$ ($\dim X$)-dimensional $1$-handlebodies smoothly embedded in $X$ (resp. which are PL submanifolds of $X$) enjoying the following properties is said to {\it define a smooth} (resp. {\it PL}) {\it multisection} of $X$ of degree $k$ where $\{X_j\}_{j=1}^k$ denotes the sequence of all our handlebodies here.
			\begin{enumerate}
				\item
				\label{def:3.1.1}	
				$X={\bigcup}_{j=1}^k X_j$.
				\item
				\label{def:3.1.2}
				${\rm Int}\ X_{j_1} \bigcap {\rm Int}\ X_{j_2}$ is empty for any pair $(j_1,j_2)$ satisfying $j_1 \neq j_2$.
				\item
				\label{def:3.1.3}
				There exists a suitable integer $l({\{X_j\}}_{j=1}^k)$ satisfying the condition $1 \leq l({\{X_j\}}_{j=1}^k) \leq k$ and $d_{l(\{X_j\}_{j=1}^k)}(k^{\prime }):=\min\{l(\{X_j\}_{j=1}^k),k^{\prime },\dim X+1\}$ is defined for each $1 \leq k^{\prime} \leq k$.  
				\item 
				\label{def:3.1.4}
				Choose distinct $1 \leq k^{\prime}<k$ numbers from $k$ integers greater than or equal to $1$ and smaller than or equal to $k$. For each integer $j^{\prime}$ of these $k^{\prime}$, we consider $X_{j^{\prime}}$ and consider the intersection. 
				Then for the previously chosen $k^{\prime}$ this is a ($\dim X-d_{l(\{X_j\}_{j=1}^k)}(k^{\prime})+1$)-dimensional compact manifold smoothly embedded in (resp. which is a PL submanifold of) $X$.
				\item 
					\label{def:3.1.5}
				Furthermore, the following conditions are hold for the submanifold of (\ref{def:3.1.4}) where $1 \leq k^{\prime \prime}<l(\{X_j\}_{j=1}^k)$
			 and $l(\{X_j\}_{j=1}^k) \leq
	k^{\prime \prime \prime}<k$ are arbitrary integers satisfying the inequalities. 
	
	The condition (\ref{def:3.1.5.3}) is for (\ref{def:3.2}) and we can also apply the presentation here.
				\begin{enumerate}
					\item \label{def:3.1.5.1}
					 Let $D_{l(\{X_j\}_{j=1}^k)}(k^{\prime \prime}):=\dim d_{l(\{X_j\}_{j=1}^k)}(k^{\prime \prime})$ in the case $(\dim X,d_{l(\{X_j\}_{j=1}^k)}(k^{\prime \prime})) \neq (2(k-1),k-1)$.
					The relation  $D_{l(\{X_j\}_{j=1}^k)}(k^{\prime \prime})<\dim X-d_{l(\{X_j\}_{j=1}^k)}(k^{\prime \prime})+1$ holds.
					In this case the submanifold obtained as the intersection of $k^{\prime \prime}$ distinct submanifolds in $\{X_j\}_{j=1}^k$ always collapses to a ($D_{l(\{X_j\}_{j=1}^k)}(k^{\prime \prime})$)-dimensional connected subpolyhedron.
					\item \label{def:3.1.5.2}
					Let $D_{l(\{X_j\}_{j=1}^k)}(k^{\prime \prime}):=\dim d_{l(\{X_j\}_{j=1}^k)}(k^{\prime \prime})-1=k-2$ in the case $(\dim X,d_{l(\{X_j\}_{j=1}^k)}(k^{\prime \prime}))=(2(k-1),k-1)$. The relation  $D_{l(\{X_j\}_{j=1}^k)}(k^{\prime \prime})<\dim X-d_{l(\{X_j\}_{j=1}^k)}(k^{\prime \prime})+1$ holds as in the previous case. In this case the submanifold obtained as the intersection of $k^{\prime \prime}$ distinct submanifolds in $\{X_j\}_{j=1}^k$ always collapses to a ($D_{l(\{X_j\}_{j=1}^k)}(k^{\prime \prime})$)-dimensional connected subpolyhedron.
					
					\item \label{def:3.1.5.3}
					 The intersection of $k^{\prime \prime}$ distinct  submanifolds in $\{X_j\}_{j=1}^k$ is always a non-empty subset of $X$ and a ($\dim X-l(\{X_j\}_{j=1}^k)+1$)-dimensional smooth (resp. PL) closed submanifold of $X$ and if its dimension is greater than $0$ then it is connected. Furthermore, its interior is embedded in the interior ${\rm Int}\ X$ and the boundary is embedded in the boundary $\partial X$.
					\end{enumerate}

			\end{enumerate}
			Furthermore, if we can take $l(\{X_j\}_{j=1}^k):=k \leq \dim X+1$, then the smooth (resp. PL) multisection is said to be {\it normal}.
			\item
			\label{def:3.2}
			For a smooth (PL) compact and connected manifold $X$ whose boundary is not empty, a family of $k$ ($\dim X$)-dimensional $1$-handlebodies smoothly embedded in $X$ (resp. which are PL submanifolds of $X$) enjoying the following two is said to {\it define a smooth} (resp. {\it PL}) {\it multisection} of $X$ of degree $k$.
			\begin{enumerate}
				\item (\ref{def:3.1.1})--(\ref{def:3.1.5}) in the previous case are enjoyed where the notation is abused.
				\item For each connected component ${\partial}_i X$ of the boudnary $\partial X$ of $X$,
				consider $X_j \bigcap {\partial}_i X$ for each $X_j$. Then the family of the resulting $k$ manifolds defines a smooth (resp. PL) multisection of ${\partial}_i X$. 
			\end{enumerate}
			We can define a smooth (PL) {\it normal} multisection similarly.
			\item \label{def:3.3}
			For a smooth (PL), closed and connected manifold $X$, a family of $k$ ($\dim X$)-dimensional compact and connected submanifolds smoothly embedded in $X$ (resp. which are PL submanifolds of $X$) enjoying the following properties is said to {\it define a smooth} (resp. {\it PL}) {\it generalized multisection} of $X$ of degree $k$ where $\{X_j\}_{j=1}^k$ denotes the sequence of all our compact and connected manifolds here.
			
			\begin{enumerate}
				\item
				\label{def:3.3.1}	
				$X={\bigcup}_{j=1}^k X_j$.
				\item
				\label{def:3.3.2}
				${\rm Int}\ X_{j_1} \bigcap {\rm Int}\ X_{j_2}$ is empty for any pair $(j_1,j_2)$ satisfying $j_1 \neq j_2$.
			\end{enumerate}
		\item \label{def:3.4}
		 Furthermore, a smooth (PL) generalized multisection of (\ref{def:3.3}) enjoying the following properties is said to be {\it generic}.  
		\begin{enumerate}
			\item \label{def:3.4.1}
		For a suitable integer $1 \leq l(\{X_j\}_{j=1}^k) \leq k$, $d_{l(\{X_j\}_{j=1}^k)}(k^{\prime}):=\min\{l(\{X_j\}_{j=1}^k),k^{\prime},\dim X+1\}$ is defined for each $1<k^{\prime} \leq k$.  
		\item 
		\label{def:3.4.2}
		Choose distinct $1 \leq k^{\prime \prime}<k$ numbers from $k$ integers greater than or equal to $1$ and smaller than or equal to $k$. For each integer $j^{\prime}$ of these $k^{\prime \prime}$, we consider $X_{j^{\prime}}$ and consider the intersection. 
		Then for the previously chosen $k^{\prime}$ this is a ($\dim X-d_{l(\{X_j\}_{j=1}^k)}(k^{\prime})+1$)-dimensional compact manifold smoothly embedded in (resp. which is a PL submanifold of) $X$.
		\end{enumerate}
	Furthermore, suppose also that for any integer $l(\{X_j\}_{j=1}^k) \leq k^{\prime \prime \prime} \leq k$, the intersection of $k^{\prime \prime}$ distinct submanifolds in $\{X_j\}_{j=1}^k$ is always a non-empty subset of $X$ and a ($\dim X-l(\{X_j\}_{j=1}^k)+1$)-dimensional smooth (resp. PL) closed submanifold of $X$ with no boundary.
	In this special case, the smooth (resp. PL) generic generalized multisection is said to be {\it properly generic}.
			\item \label{def:3.5}
			Furthermore, for a smooth or PL properly generic generalized multisection of (\ref{def:3.4}), if we can take $l(\{X_j\}_{j=1}^k):=k \leq \dim X+1$, then the smooth or PL generalized multisection is said to be {\it normal}. 
		\end{enumerate}
		
	\end{Def}
	
	Smooth (PL) multisections are smooth (resp. PL) properly generic generalized multisections. Smooth (PL) normal multisections are smooth (resp. PL) normal generalized multisections. These follow from our definition.
	
	We can define smooth (PL) generalized multisections for cases where the manifolds have non-empty boundaries. However we do not discuss them.
	
	We can define a natural equivalence relation on the class of smooth (PL) generalized multisections by diffeomorphisms (resp. piesewise smooth homeomorphisms). More precisely, each diffeomorphsim or piesewise smooth homeomorphism maps each submanifold in the family defining some multisection to one in another submanifold in the family defining another multisection. Under the equivalence relation, we can define the notion that two smooth (resp. PL) multisection are defined to be {\it equivalent}.
	\begin{Ex}
		\label{ex:2}
		\begin{enumerate}
			\item \label{ex:2.1}
		A copy of the unit disk $D^l$ of dimension $l>0$ always admits a smooth and PL normal multisection of degree $k>0$ for any positive integer $k \leq l+1$ defined by $k$ copies of the unit disk $D^l$ in the given copy of the disk $D^l$. We present an explicit example for this.
		
		For a non-negative integer $k^{\prime}$ satisfying $0 \leq k^{\prime} \leq k$, let ${{\mathbb{R}}^l}_{k^{\prime},\leq}:=\{(x_1,\cdots,x_l) \in {\mathbb{R}}^l \mid$ For $1 \leq j \leq k^{\prime}$, $x_j \geq 0$ and $x_{k^{\prime}+1} \leq 0\}$ and ${{\mathbb{R}}^l}_{k^{\prime}}:=\{(x_1,\cdots,x_l) \in {\mathbb{R}}^l \mid$ For $1 \leq j \leq k^{\prime}$, $x_j \geq 0\}$. 
		$\{D^l \bigcap {{\mathbb{R}}^l}_{j^{\prime}-1,\leq}\}_{j^{\prime}=1}^{k-1} \sqcup \{D^l \bigcap {{\mathbb{R}}^l}_{k-1}\}$ defines a desired smooth and PL normal multisection.
		
		\item \label{ex:2.2}
		
		For example, it is a fundamental fact that for an arbitrary smooth (PL) normal multisection of a compact and connected manifold $X$ of dimension $\dim X>0$ of degree $k>0$ with $0<k \leq \dim X +1$ of a suitable class, we can take a suitable smoothly embedded copy $X_{D^{\dim X}}$ of the unit disk $D^{\dim X}$ in the inteirior ${\rm Int}\ X$ enjoying the following properties.
		\begin{enumerate}
			\item $X_{D^{\dim X}}$ has a smooth and PL normal multisection of degree $k$.
			\item The previous smooth and PL multisection of $X_{D^{\dim X}}$ is defined by the family of submanifolds each of which is a copy of the unit disk $D^{\dim X}$ being also the smooth closed submanifold of ${\rm Int}\ X$ represented as the intersection of a submanifold in the family defining the given multisection of $X$ and $X_{D^{\dim X}}$. Furthermore, the smooth and PL multisection of $X_{D^{\dim X}}$ and some smooth and PL multisection in (\ref{ex:2.1}) are smooth and PL equivalent if the dimensions of the disks and the degrees are same. 
			\end{enumerate}
		  
			\item 
		\label{ex:2.3}
			The previous two smooth and PL normal multisections give smooth and PL normal multisections of the standard spheres of the boundaries. They are defined by families of copies of unit disks. The intersections as in Definition \ref{def:3} (\ref{def:3.1.4}) are standard spheres or two-point sets (with the discrete topologies).
			
	\end{enumerate}  
	\end{Ex}


\begin{Def}
\label{def:4}
We consider a smooth (PL) multisection of a compact and connected manifold $X$ in Example \ref{ex:2} (\ref{ex:2.2}) such that we can do an iteration of $i>0$ times the procedure of removing a copy of a unit disk playing roles same as ones $X_{D^{\dim X}}$ plays there and have a new and natural smooth (resp. PL) multisection again for an arbitrary positive integer $i>0$. 
We say that the given smooth (PL) multisection is {\it neat}.
Later, we consider a slightly general smooth or PL generalized multisection of a manifold $X$ in Example \ref{ex:3} and we say that such a smooth or PL multisection is {\it near neat}.
\end{Def}
By weakening the conditions naturally and arguing suitably in Example \ref{ex:2}, we have the following example.
	\begin{Ex}
	\label{ex:3}
	
For an arbitrary smooth (PL) multisection of a compact and connected manifold $X$ of dimension $\dim X>0$ of degree $k>0$ of a suitable class, we have a suitable smoothly embedded copy $X_{D^{\dim X}}$ of the unit disk $D^{\dim X}$ in the interior ${\rm Int}\ X$ enjoying the following properties.
	\begin{enumerate}
		\item $X_{D^{\dim X}}$ has a smooth and PL properly generic multisection of degree $k$.
		\item The previous smooth and PL multisection is defined by the family of submanifolds each of which is a copy of the unit disk $D^{\dim X}$ being also the smooth closed submanifold of ${\rm Int}\ X$ represented as the intersection of a submanifold in the family defining the given smooth (resp. PL) multisection of $X$ and $X_{D^{\dim X}}$.
	\end{enumerate}

In Example \ref{ex:2} (Example \ref{ex:2} (\ref{ex:2.1})), by admitting the condition $k>l+1$ or dropping the condition that a smooth or PL multisection is normal, we also have a smooth and PL properly generic multisection which is not normal in a natural way.
We call this a {\it disk with a smooth and PL almost standard generalized multisection}.

By considering a suitable class of smooth (PL) multisections of compact and connected manifolds, we can do so that the smooth and PL multisection of the unit disk we have obtained at the beginning here and a suitable smooth and PL multisection of a disk of Example \ref{ex:2} or some disk with a smooth and PL almost standard generalized multisection
are equivalent if the dimensions of the disks are same and the degrees are also same.

We also have a smooth and PL properly generic multisection which may not be normal as in Example \ref{ex:2} (\ref{ex:2.3}).
We can also do so that the intersections as in Definition \ref{def:3} (\ref{def:3.1.4}) are standard spheres or two-point sets (with the discrete topologies) for this smooth and PL properly generic multisection of the standard sphere.
\end{Ex}

Related to Example \ref{ex:3}, we can similarly define a class as we have done respecting Example \ref{ex:2}. This is presented already in Definition \ref{def:4} as smooth or PL {\it near neat} multisections.

\begin{Prop}
	A smooth {\rm (}PL{\rm )} neat multisection is a smooth {\rm (}resp. PL{\rm )} near neat multisection.
\end{Prop}

Hereafter, we mainly consider neat ones essentially.
We omit several rigorous expositions in Example \ref{ex:3} due to this. Related to this, we also remark that Example \ref{ex:3} respects some of \cite{ogawa2} for example. This studies smooth and PL near neat multisections of $3$-dimensional closed, connected and orientable manifolds which may not be neat. Rigorous expositions are left to readers. 
	\begin{Thm}
		\label{thm:2}
		A $2$-dimensional closed, connected and orientable manifold and $3$-dimensional one admit smooth and PL neat multisections of degree $2$. The $3$-dimensional case presents a {\it Heegaard splitting}. 
	\end{Thm}
	
	Note that manifolds whose dimensions are at most $3$ are uniquely regarded as PL manifolds and smooth manifolds. This is so-called Hauptvermutung, discussed in \cite{moise}, for example. For topological theory of $3$-dimensionanl manifolds, see \cite{hempel} for example. 
	
	\begin{Thm}[\cite{rubinsteintillmann1, rubinsteintillmann2}]
		\label{thm:3}
		A {\rm (}$\dim X${\rm )}-dimensional PL, closed, connected and orientable manifold $X$ admits a PL neat multisection of degree $l=\frac{\dim X}{2}$+1 if $\dim X>0$ is even and of degree $l=\frac{\dim X +1}{2}$ if $\dim X>0$ is odd. 
	\end{Thm}

	\begin{Thm}[The $4$-dimensional case is due to \cite{gaykirby} and the $5$-dimensional case is due to \cite{lambertcolemiller}]
		\label{thm:4}
		
		A $4$-dimensional smooth, compact, connected and orientable manifold admits a smooth and PL neat multisection of degree $3$, which is a so-called {\rm trisection}. This holds for $5$-dimensional cases. 
	\end{Thm}
	
	For related theory, see also \cite{martelli} for example.
	This is in the references of \cite{rubinsteintillmann1} for example. This is on so-called {\it complexities} of PL manifolds and arguments are important in the theory of Rubinstein and Tillmann.
	
	We concentrate on the following two problems and concentrate on explicit classes and examples of smooth or PL (normal generalized) multisections starting from defining such classes.
	
	\begin{Prob}
		\label{prob:1}
		In Theorem \ref{thm:3}, can we replace $l$ by another integer under a suitable situation? We may change the conditions on the class of smooth or PL (generalized) multisections in suitable ways and discuss our problems.
	\end{Prob}
	\begin{Prob}
		\label{prob:2}
		In \cite{rubinsteintillmann1, rubinsteintillmann2}, methods of construction of PL multisections (most of which are standard) are given in a general manner. 
		They introduce generalized multisections first, methods and arguments related to which are respected in Definition \ref{def:3} of our paper for example.
		It is also presented that closed and connected PL manifolds admitting nice symmetries admit PL (generalized) multisections due to the symmetries. However, these arguments do not mean presenting explicit examples and such explicit construction seems to be difficult in general. Can we give nice examples systematically starting from defining nice explicit classes?
	\end{Prob}
	The following class shows one of new work in our paper.
	
	Hereafter, we need some notions and methods from elementary algebraic topology. See \cite{hatcher} for them. 
	We introduce a class. This is an important ingredient of our new work.

		\begin{Def}
		\label{def:5}
	
			For a smooth (PL), closed and connected manifold $X$, suppose that a family $\{X_j\}_{j=1}^k$ of $k>0$ ($\dim X$)-dimensional compact and connected submanifolds smoothly embedded in $X$ (resp. which are PL submanifolds of $X$) exists.
			
			Suppose also that there exist a smooth (resp. PL), closed and connected manifold $Y$ and a smooth (resp. PL) multisection of $Y$ defined by a family $\{Y_j\}_{j=1}^k$ of $1$-handlebodies embedded in $Y$.
			 
			 The family $\{X_j\}_{j=1}^k$ is said to {\it define a smooth} (resp. {\it PL}) {\it near multisection} of $X$ of degree $k$
			 if it enjoys the following properties where we abuse the notation in Definition \ref{def:3}	(\ref{def:3.1}) for example.  
				\begin{enumerate}
					\item 	\label{def:5.0}	
					$\{X_j\}_{j=1}^k$ defines a smooth (resp. PL) generalized multisection.
				\item
				\label{def:5.1}	
				$\dim Y \leq \dim X$.
				\item
				\label{def:5.2}
					Choose distinct $1 \leq k^{\prime}<k$ numbers from $k$ integers greater than or equal to $1$ and smaller than or equal to $k$. For each integer $j^{\prime}$ of these $k^{\prime}$, we consider $X_{j^{\prime}}$ and consider the intersection. 
				Then for the previously chosen $k^{\prime}$ this is a ($\dim X-d_{l(\{Y_j\}_{j=1}^k)}(k^{\prime})+1$)-dimensional compact manifold smoothly embedded in (resp. which is a PL submanifold of) $X$.
				 Moreover, for this submanifold, the following hold where $1 \leq k^{\prime \prime}<l(\{Y_j\}_{j=1}^k)$
				 and $l(\{Y_j\}_{j=1}^k) \leq
				 k^{\prime \prime \prime}<k$ are arbitrary integers satisfying the inequalities. 
				\begin{enumerate}
					\item \label{def:5.2.1}
					
					The relation  $D_{l(\{Y_j\}_{j=1}^k)}(k^{\prime \prime})<\dim Y-d_{l(\{Y_j\}_{j=1}^k)}(k^{\prime \prime})+1$ holds.
					The submanifold obtained as the intersection of $k^{\prime \prime}$ distinct submanifolds in $\{X_j\}_{j=1}^k$ always collapses to a ($D_{l(\{Y_j\}_{j=1}^k)}(k^{\prime \prime})+\dim X-\dim Y$)-dimensional connected subpolyhedron.

					\item \label{def:5.2.2}
					The intersection of $k^{\prime \prime \prime}$ distinct submanifolds in $\{X_j\}_{j=1}^k$ is always a non-empty subset of $X$ and a ($\dim X-l(\{X_j\}_{j=1}^k)+1$)-dimensional smooth (resp. PL) closed submanifold of $X$ with no boundary and if its dimension is greater than $0$, then it is connected. 
				\end{enumerate}
			
				\item
				\label{def:5.3}
					Let $J_0$ denote the set of all integers greater than or equal to $1$ and smaller than or equal to $k$.
			Let $J$ be an arbitrary subset of the set $J_0$. Put $X_J:={\bigcap}_{j \in J} X_j$ and $Y_J:={\bigcap}_{j \in J} Y_j$. The following properties are enjoyed.
			\begin{enumerate}
				\item \label{def:5.3.1}
			 For an arbitrary subset $J$ of the set $J_0$, the fundamental groups of $X_J$ and $Y_J$ are always isomorphic.

			 	\item  \label{def:5.3.2}
			 	Let ${J_1}$ denote an arbitrary subset of the set $J_0$ whose size is smaller than $l(\{Y_j\}_{j=1}^k)$.
			 The homology group $H_{i_{J_{1,1}}}(X_{J_1};\mathbb{Z})$ is not the trivial group if and only if $i_{J_{1,1}}$ is equal to $i_{J_{1,2}}$ or $i_{J_{1,2}}+\dim X-\dim Y$ for some integer $i_{J_{1,2}} \neq 0$ such that the homology group $H_{i_{J_{1,2}}}(Y_{J_1};\mathbb{Z})$ is not the trivial group. Furthermore, for such $i_{J_{1,2}}$, the homology groups $H_{i_{J_{1,2}}}(X_{J_1};\mathbb{Z})$ and $H_{i_{J_{1,2}}+\dim X-\dim Y}(X_{J_1};\mathbb{Z})$ are isomorphic to $H_{i_{J_{1,2}}}(Y_{J_1};\mathbb{Z})$.
			 \item  \label{def:5.3.3}
			 Let ${J_2}$ denote an arbitrary subset of the set $J_0$ whose size is greater than or equal to
			  $l(\{Y_j\}_{j=1}^k)$.
			The homology group $H_{i_{J_{2,1}}}(X_{J_2};\mathbb{Z})$ is not the trivial group if and only if $i_{J_{2,1}}$ is equal to $i_{J_{2,2}}$ or $i_{J_{2,2}}+\dim X-\dim Y$ for some integer $i_{J_{2,2}} \neq 0$ such that the homology group $H_{i_{J_{2,2}}}(Y_{J_2};\mathbb{Z})$ is not the trivial group or $i_{J_{2,1}}=\dim X_{J_0}$. Furthermore, for such $i_{J_{2,2}}$, the homology groups $H_{i_{J_{2,2}}}(X_{J_2};\mathbb{Z})$ and $H_{i_{J_{2,2}}+\dim X-\dim Y}(X_{J_2};\mathbb{Z})$ are isomorphic to $H_{i_{J_{2,2}}}(Y_{J_2};\mathbb{Z})$. In addition, $H_{\dim X_{j_2}}(X_{J_2};\mathbb{Z})$ is isomorphic to $\mathbb{Z}$ if  $X_{J_2}$ is orientable. Furthermore, $H_{\dim X_{j_2}}(X_{J_2};\mathbb{Z})$ is the trivial group if $X_{J_2}$ is not orientable.
			 	\end{enumerate}

			\end{enumerate}
			
	\end{Def}

    We explain about a smooth or PL near multisection in a specific case.
    By our definitions, we have the following easily for example.
    \begin{Cor}
    	\label{cor:1}
    	\begin{enumerate}
    		\item 
    	    A smooth {\rm (}PL{\rm )} multisection is by its definition a smooth  {\rm (}resp. PL{\rm )} near multisection.	
    \item A smooth {\rm (}PL{\rm )} near multisection is by its definition a smooth {\rm (}resp. PL{\rm )} properly generic generalized multisection.
    \end{enumerate}
    \end{Cor} 
As another fact, a smooth or PL near multisection is, by the definition, not a smooth or PL multisection in general in the case $\dim Y<\dim X$. In some specific cases such as the case of Corollary \ref{cor:2}, it is a smooth or PL multisection. 

We discuss further. For example, by the definitions of a smooth or PL
    multisection and a smooth or PL near multisection and the conditions on the fundamental
    groups and the homology groups of each manifold $Y_j$,
     $Y_j$ collapses to a polyhedron whose dimension is smaller than $\dim Y_j$ by $d \geq 2$ if $\dim Y \geq 3$. This means that in the case $\dim Y \geq 3$ and the given smooth or PL multisection of $Y$ is normal, $\{X_j\}_{j=1}^k$ defines a smooth or PL generalized multisection for Definition 6 of \cite{rubinsteintillmann1}.
   In this case, this family $\{X_j\}_{j=1}^k$ also defines a smooth or PL normal generalized multisection of Definition \ref{def:3}.
	
	\section{On Main Theorems.}
    We prove Main Theorems.
	
	\begin{proof}[A proof of Main Theorem \ref{mthm:1}]
		Proposition \ref{prop:1} is regarded as a theorem on the structure of a special generic map $f:M \rightarrow {\mathbb{R}}^n$ and abuse the notation.
		
	    $Y:=W_f$ is a manifold admitting a smooth or PL multisection of degree $k>0$ such that the smooth or PL multisection induced on each connected component of $\partial Y$, which is a standard sphere, is as in Example \ref{ex:2} (\ref{ex:2.3}) or Example \ref{ex:3} from the assumption that the smooth or PL multisection is near neat.
	    
	    Investigate $X_j:={q_f}^{-1}(Y_j)$ and the intersections of finitely many subsets in the family $\{X_j\}_{j=1}^k$.
	    
	    Let $J$ be an arbitrary subset of the set $J_0$ of all integers greater than or equal to $1$ and smaller than or equal to $k$. Let $|J|$ denote the size. Put $X_J:={\bigcap}_{j \in J} X_j$ and $Y_J:={\bigcap}_{j \in J} Y_j$. We have ${q_f}^{-1}(Y_J)=X_J$ and $\dim X_J=\dim Y_J+\dim X-\dim Y$.
	    
	    $X:=M$ is the union of all subsets in the family $\{X_j\}_{j=1}^k$.

	     We will see that this defines a desired smooth or PL near multisection.

 From Proposition \ref{prop:1} with some fundamental properties of arguments on differential topology, $X_J$ can be regarded as a smooth or PL, compact, connected and orientable manifold. 
 If the smooth (PL) multisection of $Y$ is given, then $X_J$ is smooth (resp. PL).
 We can easily see that $\{X_j\}_{j=1}^k$ defines a smooth (resp. PL) generalized multisection.
 We have the property (\ref{def:5.0}) of Definition \ref{def:5}. We can also see that the property (\ref{def:5.1}) of Definition \ref{def:5} is enjoyed. 
	    
	    The property (\ref{def:5.3}) of Definition \ref{def:5} are all shown via elementary algebraic topological arguments as follows.
	    
	   $X_J$ is regarded as a smooth or PL manifold and compact, connected and orientable. 
	    $M$ and $W_f$ are orientable and each $Y_j$ in the family $\{Y_j\}_{j=1}^k$ is a 1-handlebody.
	    
	    As in (\ref{def:5.3.3}) here, we consider the case $J_2 \subset J_0$.

	    $X_{J_2}$ is, as a closed PL manifold, PL homeomorphic to one obtained by attaching $\partial Y_{J_2} \times D^{m-n+1}$ to $Y_{J_2} \times S^{m-n}$ by the product map of two PL homeomorphisms between the boundaries. We also have the property (\ref{def:5.2.2}) of Definition \ref{def:5}.
	    
	    In the case $J_1 \subset J_0$ as in (\ref{def:5.3.2}) here,
	    $X_{J_1}$ is, as a PL manifold, PL homeomorphic to one attaching a copy of the unit disk $D^{\dim \partial Y_{J_1}+m-n}$ along $D_{Y_{J_1},\dim \partial Y_{J_1}} \times S^{m-n} \subset Y_{J_1} \times S^{m-n}$ where $D_{Y_{J_1},\dim \partial Y_{J_1}}$ is a copy of the ($\dim \partial Y_{J_1}$)-dimensional unit disk $D^{\dim \partial Y_{J_1}}$ embedded smoothly (resp. as a PL submanifold) in $\partial Y_{J_1}$.

	    Again, the assumption that $Y_j$ is a 1-handlebody and the orientability of $X_j$ are essential in discussing the structures of $X_j$ and $X_J$ here for any subset $J \subset J_0$. 
	    We need to apply K\"unneth theorem for the product $Y_J \times S^{m-n}$. We also need Mayer-Vietoris sequence and Seifert-van Kampen theorem for the pair of the copy of the unit disk $D^{\dim \partial Y_{J_1}+m-n}$ and $Y_{J_1} \times S^{m-n}$, glued in the represented way to present $X_{J_1}$ in the case $J_1 \subset J_0$. In the case $J_2 \subset J_0$, we consider similar arguments for the pair of a manifold diffeomorphic to the disjoint union of two copies of the $D^{m-n}$ or the product of a standard sphere and $D^{m-n}$ and $Y_{J_2} \times S^{m-n}$, glued in the represented way to present a closed manifold $X_{J_2}$. Remember that the smooth and PL multiesctions as in Example \ref{ex:2} (\ref{ex:2.3}) or Example \ref{ex:3} are induced on each connected component of the boundary $\partial Y_{J_2}$, which is assumed to be a standard sphere.
	    
	     We see that our properties (on homology groups and fundamental groups in) (\ref{def:5.3}) of Definition \ref{def:5} are enjoyed.
	    
	    We concentrate on the case $J_1 \subset J_0$ again and show the property (\ref{def:5.2}) of Definition \ref{def:5}.
	    By the attachment of a copy of the disk to $Y_{J_1} \times S^{m-n}$ before, $X_{J_1}$ collapses to a polyhedron whose dimension is at least $\max\{ D_{l(\{Y_j\}_{j=1}^k)}(|J_1|)+m-n,1+m-n\}= D_{l(\{Y_j\}_{j=1}^k)}(|J_1|)+m-n$ since $ D_{l(\{Y_j\}_{j=1}^k)}(|J_1|)>0$ is known by the assumption. 
	   More precisely, this collapses to a connected polyhedron PL homoemorphic to one obtained by attaching a copy of $D^1 \times S^{m-n}$ and a copy of the unit disk $D^{m-n+1}$ to the product of the ($D_{l(\{Y_j\}_{j=1}^k)}(|J_1|)$)-dimensional compact and connected polyhedron $Y_{J_1}$ is assumed to collapse to and $S^{m-n}$. We have the property (\ref{def:5.2.1}). The property (\ref{def:5.2.2}) is already shown. We also have the property (\ref{def:5.2}).
	   
	    We can prove the smooth case and the PL case in this way. We can also prove the case where the smooth or PL multisection is neat in a similar way. This completes the proof.

	\end{proof}
The construction in our proof yields the following corollary.

\begin{Cor}
	\label{cor:2}
	If ${\{Y_j\}}_{j=1}^k$ in Main Theorem \ref{mthm:1} defines a smooth and PL generalized multisection which is equivalent to a smooth and PL normal multisection of the unit disk $D^l$ of Example \ref{ex:2} {\rm (}\ref{ex:2.1}{\rm )}, then the resulting smooth and PL near multisection is also normal and equivalent to one in Example \ref{ex:2} {\rm (}\ref{ex:2.1}{\rm )}.
\end{Cor}

\begin{MainThm}
	\label{mthm:2}
	Let $m>n \geq 1$ be integers.
	\begin{enumerate}
		\item \label{mthm:2.0} Let $n=1$. If an $m$-dimensional closed, connected and orientable manifold $M$ admits a special generic map $f:M \rightarrow \mathbb{R}$, then it admits a PL and smooth normal near multisection of degree $1$.
		\item \label{mthm:2.1}
		 If an $m$-dimensional closed, connected and orientable manifold $M$ admits a special generic map $f:M \rightarrow {\mathbb{R}}^n$ such that in Proposition \ref{prop:1}, the boundary $\partial W_f \subset W_f$ is a disjoint union of standard spheres, then it admits a PL normal near multisection of degree $\frac{n}{2}+1$ if $n$ is even and of degree $\frac{n+1}{2}$ if $n$ is odd.
		\item \label{mthm:2.2}
		 Let $n=2$. If an $m$-dimensional closed, connected and orientable manifold $M$ admits a special generic map $f:M \rightarrow {\mathbb{R}}^n$, then it admits a PL and smooth normal near multisection of degree $2$.
		\item \label{mthm:2.3}
		Let $n=3$. If an $m$-dimensional closed, connected and orientable manifold $M$ admits a special generic map $f:M \rightarrow {\mathbb{R}}^n$ such that in Proposition \ref{prop:1}, the boundary $\partial W_f \subset W_f$ is a disjoint union of $2$-dimensional standard spheres, then it admits a PL and smooth normal near multisection of degree $2$.
		\item \label{mthm:2.4}
		 Let $n=4,5$. If an $m$-dimensional closed, connected and orientable manifold $M$ admits a special generic map $f:M \rightarrow {\mathbb{R}}^n$ such that in Proposition \ref{prop:1}, the boundary $\partial W_f \subset W_f$ is a disjoint union of {\rm (}$n-1${\rm )}-dimensional standard spheres, then it admits a PL and smooth normal near multisection of degree $3$.
	\end{enumerate}
\end{MainThm}
\begin{proof}[A proof of Main Theorem \ref{mthm:2}]
	 The case (\ref{mthm:2.0}) is easily shown by considering Morse functions of the Reeb's theorem.
	We can show the remaining four as an application of Main Theorem \ref{mthm:1} with some arguments and known results from special generic maps.
	 The case (\ref{mthm:2.1}) is due to Theorem \ref{thm:3}. The case (\ref{mthm:2.2}) is closely related to Theorem \ref{thm:1} (\ref{thm:1.1}) and due to Theorem \ref{thm:2}. We have this with the easily known fact that we have a $2$-dimensional case of Theorem \ref{thm:2} and the orientability of $W_f$ in Proposition \ref{prop:1} for example. The orientablity is due to the fact that a manifold smoothly immersed into the Euclidean space of the same dimension must be orientable. The case (\ref{mthm:2.3}) is also due to Theorem \ref{thm:2}. 
	Note again that the $2$-dimensional and $3$-dimensional cases can be discussed in the smooth category essentially by Hauptvermutung, discussed in \cite{moise} for example. The case (\ref{mthm:2.4}) is due to Theorem \ref{thm:4} and the orientability of $W_f$ in Proposition \ref{prop:1}.
	
	This completes the proof.
\end{proof}

Main Theorem \ref{mthm:2} presents PL or smooth near multisections explicitly. As an answer to Problem \ref{prob:1}, we can construct one which can not be a smooth or PL multisection such that $l$ there is lower than the integer easily. Main Theorem \ref{mthm:2} also gives an answer to Problem \ref{prob:2}. Example \ref{ex:4} explains about these cases.

\begin{Ex}
	\label{ex:4} 
	Special generic maps in Theorem \ref{thm:1} are all for Main Theorem \ref{mthm:2} (\ref{mthm:2.2}) and (\ref{mthm:2.3}).
	
    We can consider a case such that $W_f$ in Proposition \ref{prop:1} is a $4$-dimensional compact and simply-connected manifold which is obtained by removing the interiors of finitely many smoothly and disjointly embedded copies of the unit disk $D^4$ from a smooth manifold represented as a connected sum of $k$ copies of $S^2 \times S^2$ for any positive integer $k \geq 1$: the connected sum is considered in the smooth category. Consider the case where the manifold of the domain is $5$-dimensional. The resulting manifold of the domain of the special generic map is, according to \cite{nishioka}, diffeomorphic to a manifold represented as a connected sum of $2k$ manifolds each of which is the total space of a linear bundle over $S^2$ whose fiber is the $3$-dimensional unit sphere $S^3$. This is for some of Main Theorem \ref{mthm:2} (\ref{mthm:2.4}). Note also that there exist exactly two types of linear bundles over $S^2$ whose fibers are the $3$-dimensional unit sphere $S^3$ and that the types of the linear bundles here can be arbitrary.
    
    Let $G_1$ be a finitely generated free commutative group. Let $G_2$ be a finite commutative group which is represented as the direct sum of two copies of some group. 
    We can consider a case such that $W_f$ in Proposition \ref{prop:1} is a $5$-dimensional compact and simply-connected manifold which is obtained by removing the interiors of finitely many smoothly and disjointly embedded copies of the unit disk $D^5$ from a $5$-dimensional smooth closed and simply-connected manifold $X$ which is a so-called {\it spin} manifold such that the homology group $H_j(X;\mathbb{Z})$ is isomorphic to the direct sum $G_1 \oplus G_2$ for $j=2$, $G_1$ for $j=3$ and the trivial group otherwise.
    See also \cite{barden} for complete classifications of $5$-dimensional closed and simply-connected manifolds in the topology, PL and smooth categories, which are a key ingredient in \cite{nishioka}, determining $5$-dimensional closed and simply-connected manifolds admitting special generic maps into arbitrary Euclidean spaces completely.
    
    Here, we consider the case where the manifold of the domain is $6$-dimensional closed and connected one.
     In this case, the resulting manifold $M$ of the domain of the special generic map is, according to arguments and results of \cite{kitazawa5} for example, some $6$-dimensional spin, closed and simply-connected manifold which is in the classifications of some studies such as \cite{jupp,wall,wall2,zhubr1,zhubr2}. This is for some of Main Theorem \ref{mthm:2} (\ref{mthm:2.4}). Moreover, for example, for such a $6$-dimensional spin, closed and simply-connected manifold $M$, the homology group $H_2(M;\mathbb{Z})$ is isomorphic to the direct sum $G_1 \oplus G_2$ for $j=2$, the direct sum $G_1 \oplus G_1 \oplus G_2$ for $j=3$ and we can know the others by Poincar\'e duality theorem for $M$ in the case where the boundary $\partial W_f$ of the manifold $W_f$ is connected and a $4$-dimensional standard sphere.  
\end{Ex}
\begin{Rem}
	For example, in Example \ref{ex:4}, we do not know whether the $6$-dimensional closed and simply-connected manifolds admit smooth or PL multisections of degree $3$. They admit PL neat multisections of degree $4$ by \cite{rubinsteintillmann1, rubinsteintillmann2}.
\end{Rem}
Last, we present another problem on classes of smooth or PL generalized multisections.

\begin{Prob}
	In our paper, we have introduced various classes of smooth or PL generalized multisections based on existing studies. Can we present nice and meaningful examples for the classes. Caw we give better classes and study them?
\end{Prob}
	
	\end{document}